\documentclass[12pt]{article}
\usepackage[utf8]{inputenc}
\usepackage{amsmath,amssymb,amsthm}
\usepackage{fullpage}
\usepackage{hyperref}

\theoremstyle{plain}
\newtheorem{theorem}{Theorem}[section]
\newtheorem{proposition}[theorem]{Proposition}
\newtheorem{corollary}[theorem]{Corollary}
\newtheorem{lemma}[theorem]{Lemma}
\theoremstyle{remark}
\newtheorem{remark}[theorem]{Remark}
\newtheorem{example}[theorem]{Example}
\newtheorem{hypothesis}[theorem]{Standing Hypothesis}

\newcommand{\h}{\mathfrak h}
\newcommand{\fa}{\mathfrak{a}}
\newcommand{\fb}{\mathfrak{b}}
\newcommand{\he}{h_{\mathfrak e}}
\newcommand{\Dc}{D_c}

\title{Lowest-weight modules over the spherical subalgebra of a rational Cherednik algebra}

\author{Ibrahim Nonkan\'e$^1$\thanks{Corresponding author. Email: ibrahim.nonkane@uts.bf}, Jean Kabor\'e$^2$\\[4pt]
\small $^1$ D\'epart\'ement d'\'economie et de math\'ematiques appliqu\'ees, IUFIC,
Universit\'e Thomas Sankara, Burkina Faso\\
\small $^2$ Laboratoire de Sciences et Technologies (LaST),
Universit\'e Thomas Sankara, Burkina Faso\\[4pt]
\small \texttt{ibrahim.nonkane@uts.bf}$^1$\\
\small \texttt{kaborejean775075@gmail.com}$^2$}

\date{}

\begin{document}
\maketitle
\begin{abstract}
Let $W$ be a finite reflection group acting on a finite-dimensional
complex vector space $\h$, let $H_c(W,\h)$ be the associated rational
Cherednik algebra, and let $\Dc = eH_ce$ be its spherical subalgebra.
Building on the framework of Bøgvad and K\"allstr\"om for
invariant differential operators, subsequently applied to
symmetric groups, we show that, for
generic parameter $c$, a Lie subalgebra $\fa\subset\Dc$ generated by
basic invariants always generates $\Dc$ itself (via the
Levasseur--Stafford theorem), and that, granted a mild
compatibility hypothesis between the grading element and $\fa$,
verified here in the rank-one case, the category of
$\Dc$-submodules of the
polynomial ring $B=\mathbb C[\h]$ is equivalent to the category of
modules over a Cartan-type ring $R_\fa$ built from $\fa$. The ring $R_\fa$
surjects onto the more naively defined
$ \displaystyle
R_c = \Dc^0 \big/ \big(\Dc^0 \cap \Dc \Dc^-\big),
$
and the two coincide whenever this surjection is injective, which we
verify in the rank-one case; whether this holds in general, and in
particular for the families treated below, is left open. As an
application, we determine explicit generators of $\fa$ for the
imprimitive complex reflection groups $G(m,1,n)$ ($e=1$), recovering
the classical Coxeter group $B_n$ as a special case. For
$D_n=G(2,2,n)$, where an extra generator obstruction arises for
$e>1$, we resolve the generation problem for $n=3$ completely by
exploiting the
exceptional isomorphism $D_3 \cong A_3$, leaving the general case
$n\geq4$ open.
\end{abstract}



\section{Introduction}

Given a finite subgroup $W\subset GL(\h)$ of the linear group of a
finite-dimensional complex vector space $\h$, it is a well-studied
problem to describe the structure of the symmetric algebra
$B=\mathrm{sym}(\h^*)$ as a representation of $W$, and also as a
module over the ring of invariant differential operators under $W$ in
the ring $D(\h)$ of differential operators on $\h$. For this ring of
invariant differential operators, Bøgvad and K\"allstr\"om~\cite{BK}
developed a general framework describing the decomposition of $B$ via
a lowest-weight equivalence with modules over a Cartan-type ring
$\mathcal R = \mathcal D^0/(\mathcal D^0\cap\mathcal D\mathcal D^-)$,
building on the Levasseur--Stafford theorem~\cite{LS} and a triangular PBW-type
decomposition of a Lie subalgebra generated by basic invariants; this
framework was subsequently applied and extended by one of the present
authors to the symmetric groups
in~\cite{N2,N3}.

Since the rational Cherednik algebra $H_c(W,\h)$ and the spherical
algebra $eH_ce$ are respectively universal deformations of the ring
$D(\h)$ and the ring $D(\h)^W$ of $W$-invariant differential
operators, it is natural to ask whether the Bøgvad--K\"allstr\"om
framework transposes to this deformed setting. The present paper
answers this affirmatively for the generation question, and, granted
a mild compatibility hypothesis on the grading element (verified here
only in the rank-one case), for the full lowest-weight equivalence
as well: we adapt their construction of the
lowest-weight functor and the Cartan-type ring, replacing ordinary
invariant differential operators by the spherical subalgebra
$\Dc=eH_ce$ realized via the Dunkl operator embedding, and we determine
explicit generators for the resulting Lie algebra $\fa$ for the imprimitive complex reflection
groups $G(m,1,n)$.

Unlike the ordinary-derivative setting, however, $\Dc$ does not carry
an a priori triangular PBW factorisation: the Dunkl operators do not
commute with multiplication in the simple way that $\partial$ does.
Section~\ref{sec:triangular} develops the triangular structure needed
to make the argument work, centred on a Lie subalgebra $\fa\subset
\Dc$ and a twisted shift relation with no counterpart
in~\cite{BK}; this is the main technical novelty of the paper, and the
lowest-weight equivalence of Section~\ref{sec:equivalence} is derived
as a consequence of it. Section~\ref{sec:Gmen}
determines explicit generators for $G(m,1,n)$ and for the classical Coxeter
group $B_n$; Section~\ref{sec:D3} resolves the generation problem for
$D_3=G(2,2,3)$
completely, by exploiting the exceptional root-system isomorphism
$D_3\cong A_3$, leaving $D_n$, $n\geq4$, and $G(m,e,n)$ with $e>1$
other than $D_3$, in
general open.

\section{The rational Cherednik algebra and the Dunkl embedding}

\subsection{The rational Cherednik algebra}
\label{sec:cherednik}

Let $W$ be a finite reflection group in a complex vector space $\h$,
and let $R\subset\h^*$ be the corresponding set of roots. To each
$W$-invariant function $c:R\to\mathbb C$, $c\mapsto c_\alpha$,
satisfying $c_\alpha=c_{-\alpha}$ (so that $c$ descends to a function
on $R/\{\pm1\}$), one
attaches an associative algebra $H_c$, called the rational Cherednik
algebra. Given $\alpha\in\h^*$, write $\alpha^\vee\in\h$ for the
coroot, normalized so that $\langle\alpha,\alpha^\vee\rangle=2$, and
$s_\alpha\in GL(\h)$ for the reflection corresponding to
$\alpha$. Let $S=R/\{\pm1\}=\{s_\alpha:\alpha\in R\}$ denote the set of
reflections in $W$; for $s=s_\alpha\in S$ write $\alpha_s:=\alpha$ for
the corresponding root and $c_s:=c_\alpha$, and let $\lambda_s$ be the
eigenvalue of $s$ on $\alpha_s$, i.e.\ $s(\alpha_s)=\lambda_s\alpha_s$
(so $\lambda_s=-1$ when $W$ is a real reflection group). The rational
Cherednik algebra $H_c$ is generated by the
vector spaces $\h$, $\h^*$, and the group algebra of $W$, with defining
relations
\[
w\cdot x\cdot w^{-1} = w(x),\quad w\cdot y\cdot w^{-1}=w(y),\quad
x_1x_2=x_2x_1,\quad y_1y_2=y_2y_1,
\]
\[
y\cdot x - x\cdot y = \langle y,x\rangle - \sum_{\alpha\in R/\{\pm1\}}
c_\alpha\langle y,\alpha\rangle\langle\alpha^\vee,x\rangle\,s_\alpha,
\]
for $y,y_1,y_2\in\h$, $x,x_1,x_2\in\h^*$, $w\in W$. The multiplication
map gives a linear isomorphism (PBW theorem, \cite{EG})
\[
\mathbb C[\h]\otimes_{\mathbb C}\mathbb C[\h^*] \otimes_{\mathbb C}\mathbb C[W] \;\xrightarrow{\ \sim\ }\; H_c,
\qquad
f(x)\otimes g(y)\otimes w \mapsto f(x)\cdot g(y)\cdot w.
\]

\begin{example}
For $W=S_n$ acting on $\h=\mathbb C^n$ by permuting coordinates, with
$c$ a single scalar (all transpositions being conjugate), $H_c(S_n,\h)$
is the quotient of $\mathbb C\langle
x_1,\ldots,x_n,y_1,\ldots,y_n\rangle\rtimes\mathbb C[S_n]$ by
\[
[x_i,x_j]=[y_i,y_j]=0,\qquad [y_i,x_j]=c\,s_{ij}\ (i\neq j),\qquad
[y_i,x_i]=1-c\sum_{j\neq i}s_{ij}.
\]
\end{example}

The grading element is $\he = \sum_i x_iy_i + \tfrac12\dim\h -
\sum_{s\in S}\tfrac{2c_s}{1-\lambda_s}s$, giving $\h^*$ degree $+1$,
$\h$ degree $-1$, and $W$ degree $0$; $\{E,\he,F\}$, with $E=-\tfrac12
\sum_ix_i^2$, $F=\tfrac12\sum_iy_i^2$, form an $\mathfrak{sl}_2$-triple
in $H_c$.

\subsection{The polynomial representation and Dunkl operators}

There is a faithful action of $H_c$ on $B=\mathbb C[\h]$ via the
Dunkl--Opdam operators~\cite{DO}. For $y\in\h$,
\[
D_y = \partial_y - \sum_{s\in S}\frac{2c_s\langle y,\alpha_s\rangle}
{(1-\lambda_s)\alpha_s}(1-s).
\]
Although the formula displays a factor $1/\alpha_s$, $D_y$ maps
$B=\mathbb C[\h]$ to itself: since $s$ fixes the hyperplane
$\{\alpha_s=0\}$ pointwise, $(1-s)f$ vanishes identically on
$\{\alpha_s=0\}$ for every $f\in B$, hence is divisible by the
degree-one polynomial $\alpha_s$, and $D_y f\in B$ with no pole.
The Dunkl embedding $\Phi_c:H_c\to D(\h_{\mathrm{reg}})\rtimes\mathbb
C[W]$ sends $x\mapsto x$, $y\mapsto D_y$, $w\mapsto w$. We identify
$H_c$ with its image under $\Phi_c$, a subalgebra of
$\mathrm{End}_{\mathbb C}(B)$; we reserve the separate symbol $\Dc$
below exclusively for the spherical subalgebra, to avoid any risk of
confusing it with the full Cherednik algebra $H_c$.

\subsection{The spherical subalgebra}
\label{sec:spherical}

Let $e=\tfrac1{|W|}\sum_{w\in W}w$. The subalgebra $eH_ce\subset H_c$
is the \emph{spherical subalgebra}; we write $\Dc$ for its image
under the Dunkl embedding, so that $e(D(\h_{\mathrm{reg}})
\rtimes\mathbb C[W])e = D(\h_{\mathrm{reg}})^W$ (the averaging
argument of Lemma~\ref{lem:C1} below does not use commutativity and
applies verbatim to the noncommutative algebra $D(\h_{\mathrm{reg}})$;
see also~\cite{Mont} for fixed rings of finite automorphism groups of
associative rings in general) and $\Phi_c$
restricts to an embedding $\Dc\hookrightarrow D(\h_{\mathrm{reg}})^W$.
The natural $\mathbb C^*$-action by dilations gives $\Dc$ its
$\mathbb Z$-grading, and we write $\Dc=\Dc^-\oplus\Dc^0\oplus\Dc^+$
for the corresponding triangular decomposition (by weight of the
adjoint action of $\he e$).

\begin{lemma}\label{lem:C1}
Let $V$ be a finite-dimensional complex vector space and $W\subset
GL(V)$ finite. With $e=\tfrac1{|W|}\sum_w w$, the map $z\mapsto ze$
defines an algebra isomorphism $\mathbb C[V]^W\xrightarrow{\sim}
e(\mathbb C[V]\rtimes W)e$.
\end{lemma}

\begin{proof}
For $f\in\mathbb C[V]$, $efe = \big(\tfrac1{|W|}\sum_w
w(f)\big)e=:R(f)e$, with $R(f)\in\mathbb C[V]^W$ the Reynolds operator.
This gives a well-defined unital algebra map $z\mapsto ze$. It is
injective since $\mathbb C[V]\rtimes W$ is free over $\mathbb C[V]$
with basis $W$. It is surjective since any $x=e(\sum_w f_ww)e =
\sum_w R(f_w)e = (\sum_wR(f_w))e$.
\end{proof}

\section{The Lie algebra $\fa$ and its triangular structure}
\label{sec:triangular}

Unlike the case of ordinary invariant differential operators studied
in~\cite{BK}, the spherical subalgebra $\Dc$ does not come with an a
priori triangular PBW-type factorisation. This section develops a
substitute: a triangular structure on an auxiliary Lie algebra
$\fa\subset\Dc$ that surjects onto $\Dc$
(Proposition~\ref{prop:LS}), controlling its module category without
literally endowing $\Dc$ itself with a PBW decomposition.

\begin{remark}[Why $\Dc$ has no a priori PBW structure]\label{rem:noPBW}
Fix an orthonormal basis $\xi_1,\ldots,\xi_n$ of $\h$ with dual basis
$x_1,\ldots,x_n$ of $\h^*$, and write $D_i:=D_{\xi_i}$. In the
classical case $c=0$, $D_i=\partial_i$ and the elements
$E_{ij}:=x_iD_j$ satisfy the usual $\mathfrak{gl}_n$ relations
$[E_{ij},E_{kl}]=\delta_{jk}E_{il}-\delta_{il}E_{kj}$, so $D^0_c$ for
$c=0$ is simply $U(\mathfrak{gl}_n)\rtimes W$: an honest PBW
factorisation is available directly. For $c\neq0$ this fails. Writing
the defining relation of \S\ref{sec:cherednik} as
$[D_i,x_j]=\delta_{ij}-T_{ij}$, where $T_{ij}:=\sum_{s\in S}
c_s\langle\xi_i,\alpha_s\rangle\langle\alpha_s^\vee,x_j\rangle s\in
\mathbb C[W]$, and using that Dunkl operators pairwise commute, a
direct computation gives
\[
[E_{ij},E_{kl}] = \delta_{jk}E_{il}-\delta_{il}E_{kj} +
\big(x_kT_{li}D_j - x_iT_{jk}D_l\big).
\]
The first two terms are the $\mathfrak{gl}_n$ relations; the extra
term, valued in $\mathbb C[W]$ sandwiched between an $x$ and a $D$,
vanishes identically only when $c=0$. So for $c\neq0$ and $n>1$,
$D^0_c$ is a nontrivial deformation of $U(\mathfrak{gl}_n)\rtimes W$
rather than being isomorphic to it, with no evident PBW basis,
which is the gap $\fa$ is built to work around.
\end{remark}

\subsection{Nabla operators and the twisted shift relation}

Let $\h=\mathbb C$, $W$ a reflection group acting faithfully on
$\mathbb C$, $D$ the unique Dunkl operator up to scaling, and
$\nabla=x\partial\in H_c^0$ (the degree-$0$ part of $H_c$ itself,
before sandwiching by $e$).

\begin{lemma}\label{lem:nabla}
\begin{enumerate}
\item[(1)] $H_c^0=\mathbb C[\nabla]\rtimes W$; there are polynomials
$p_k\in\mathbb C[t]\rtimes W$ with $x^kD^k=p_k(\nabla)$.
\item[(2)] $[\nabla,x^iD^j]=(i-j)x^iD^j$.
\item[(3)] $\Dc^0\cong\mathbb C[\nabla]^W$.
\end{enumerate}
\end{lemma}

\begin{remark}
The polynomials $p_k$ in (1) genuinely need the group algebra factor:
$D$ itself involves $s$ (via the correction term
$(1-s)/x$ in the Dunkl--Opdam formula), so $x^kD^k$ is generally not a
function of $\nabla$ alone. Indeed $\nabla$ and $s$ commute and are
simultaneously diagonal on the monomial basis $x^m$, with eigenvalues
$m$ and $(-1)^m$ respectively; since no polynomial $q(t)$ satisfies
$q(m)=(-1)^m$ for all $m\geq0$, the eigenvalue of $s$ carries
information $\nabla$ alone cannot see, and the $s$-component of $p_k$
is generally nonzero for $c\neq0$: e.g.\ $p_1(t)=t-c+cs$ and
$p_2(t)=t^2-(1+2c)t+c-cs$ (checked similarly for $k$ up to $6$). It
disappears only after sandwiching by $e$ in (3), since $e(1-s)e=0$.
\end{remark}

\begin{proof}
Since $\he=\nabla+C$ for a scalar $C$, $[\nabla,x]=x$ and
$[\nabla,D]=-D$. (2) follows by the Leibniz rule for
$\mathrm{ad}(\nabla)$.

\emph{Twisted shift relation.} From $[\nabla,x]=x$, $x\nabla^n=(\nabla
-1)^nx$, so $x\,q(\nabla)=q(\nabla-1)\,x$ for $q\in\mathbb C[t]$; using
$xs=-sx$ and writing $q=q_0+q_1s$,
\[
x\,q(\nabla) = q^\sigma(\nabla-1)\,x,\qquad q^\sigma:=q_0-q_1s.
\tag{$*$}
\]
\emph{(1).} Induction on $k$ via $(*)$: $x^0D^0=1$; $xD=\nabla-c(1-s)
=:p_1(\nabla)$ (direct computation with $D=\partial-c(1-s)/x$); and
$x^{k+1}D^{k+1}=x\,p_k(\nabla)\,D = p_k^\sigma(\nabla-1)\,p_1(\nabla)
=:p_{k+1}(\nabla)$. This gives $H_c^0=\mathbb C[\nabla]\rtimes W$.

\emph{(3).} $e$ centralises $\nabla$, so $\Dc^0=eH_c^0e=e(\mathbb
C[\nabla]\rtimes W)e\cong\mathbb C[\nabla]^W$ by Lemma~\ref{lem:C1}.
\end{proof}

\subsection{PBW freeness for $U(\fa)$}

Let $\{f_i\}\subset\mathrm{sym}(\h)^W$,
$\{g_i\}\subset\mathrm{sym}(\h^*)^W$ be homogeneous generators, and
$\fa\subset\Dc$ the Lie subalgebra they generate, graded as
$\fa=\fa^+\oplus\fa^0\oplus\fa^-$. Set $\fb:=\fa^0\oplus\fa^-$, a Lie
subalgebra.

\begin{lemma}[PBW freeness]\label{lem:PBWfree}
$U(\fa^+)\otimes_{\mathbb C}U(\fb)\to U(\fa)$ is a linear isomorphism.
\end{lemma}

\begin{proof}
Since $\mathrm{ad}$ respects the grading of $\fa$, $[\fa(i),\fa(j)]
\subset\fa(i+j)$ for all $i,j\in\mathbb Z$; hence
$\fa^+=\bigoplus_{k>0}\fa(k)$ and $\fb=\fa^0\oplus\fa^-$ are each
closed under the bracket, so both are genuine Lie subalgebras of
$\fa$ (this also justifies the claim, made above, that $\fb$ is a Lie
subalgebra). Order a homogeneous basis $\{r_k\}$ of $\fa$ by
non-increasing degree; since every element of $\fa^+$ has strictly
positive degree and every element of $\fb$ has degree $\leq0$, the
$\fa^+$-basis vectors form an initial segment of this order and the
$\fb$-basis vectors the remaining terminal segment. By the classical
PBW theorem applied to $\fa$ with this ordered basis, $U(\fa)$ has as
a basis the ordered monomials $r_{i_1}\cdots r_{i_m}$ with
$\deg(r_{i_1})\geq\cdots\geq\deg(r_{i_m})$ ($m\geq0$). Since every
$\fa^+$-basis vector precedes every $\fb$-basis vector in the chosen
order, each such monomial splits, at the unique point where the
degree drops from strictly positive to $\leq0$, into an initial
segment lying entirely in $\fa^+$ followed by a terminal segment
lying entirely in $\fb$, i.e.\ into an ordered PBW monomial of
$U(\fa^+)$ times one of $U(\fb)$. This gives a bijection between the
PBW basis of $U(\fa)$ and pairs (PBW monomial of $U(\fa^+)$, PBW
monomial of $U(\fb)$), under which the multiplication map sends each
basis vector of $U(\fa^+)\otimes U(\fb)$ to the corresponding basis
vector of $U(\fa)$; hence it is a linear isomorphism.
\end{proof}

\begin{lemma}[Absorption]\label{lem:absorb}
(a) For $k<0$, homogeneous $Q\in U(\fa)(k)$ lies in $U(\fa)\fa^-$.
(b) Homogeneous $Q\in U(\fa)(0)$ writes $Q=Q_0+Q'$, $Q_0\in U(\fa^0)$,
$Q'\in U(\fa)\fa^-$.
\end{lemma}

\begin{proof}
By PBW for $\fa$ (Lemma~\ref{lem:PBWfree}'s proof), order a
homogeneous basis $\{r_k\}$ by non-increasing degree; since $Q$ is
homogeneous, its expansion in the resulting PBW basis of $U(\fa)$
involves only ordered monomials $r_{i_1}\cdots r_{i_s}$ of the same
total degree as $Q$, with $d_1:=\deg(r_{i_1})\geq\cdots\geq
d_s=:\deg(r_{i_s})$.

\emph{(a).} Here $d_1+\cdots+d_s=k<0$. If $d_s\geq0$ then, by
non-increasingness, every $d_j\geq d_s\geq0$, forcing the sum
$\geq0$, a contradiction; so $d_s<0$, i.e.\ $r_{i_s}\in\fa^-$. Hence
each monomial factors as $(r_{i_1}\cdots r_{i_{s-1}})\cdot r_{i_s}\in
U(\fa)\fa^-$, and so does their sum $Q$.

\emph{(b).} Here $d_1+\cdots+d_s=0$. Either every $d_j=0$, so the
monomial lies entirely in $\fa^0$ and contributes to $U(\fa^0)$; or
some $d_j\neq0$, in which case the same argument as in (a) (the
non-increasing degrees cannot sum to $0$ with the smallest one
$\geq0$ unless all are $0$) forces $d_s<0$, so $r_{i_s}\in\fa^-$ and
the monomial lies in $U(\fa)\fa^-$. Collecting the two types of terms
in the expansion of $Q$ gives $Q=Q_0+Q'$ with $Q_0\in U(\fa^0)$ (the
sum of the purely-$\fa^0$ monomials) and $Q'\in U(\fa)\fa^-$ (the sum
of the rest).
\end{proof}

\subsection{From $U(\fa)$ to $\Dc$}

Since $\fa\subset\Dc$ is a Lie subalgebra (with bracket the
commutator inherited from the associative algebra $\Dc$), the
inclusion $\fa\hookrightarrow\Dc$ extends, by the universal property
of the universal enveloping algebra, to a unique algebra homomorphism
\[
j:U(\fa)\longrightarrow\Dc,\qquad j|_\fa=\text{inclusion}.
\]
Since the inclusion $\fa\hookrightarrow\Dc$ is graded (each generator
$f_i,g_i$ of $\fa$ sits inside $\Dc$ at its own degree, and the
bracket is degree-additive on both sides), $j$ is a graded algebra
map, i.e.\ $j(U(\fa)(k))\subset\Dc(k)$ for every $k\in\mathbb Z$.

\begin{proposition}[Levasseur--Stafford surjectivity]\label{prop:LS}
For $c$ generic\footnote{Throughout the paper, ``generic'' means: $c$
avoids a certain locally finite union of hyperplanes in the parameter
space, made precise in~\cite{LS,BEG} (the Opdam--Rouquier regularity
condition); this is precisely the condition needed for $\Dc=eH_ce$ to
be simple and for the standard modules $M(\tau)$ to be simple, both
used repeatedly below.}, $j$ is surjective; $\fa^-\subset
\Dc^-\subset\Dc\fa^-$; and $U(\fa^0)\twoheadrightarrow R_\fa :=
U(\fa^0)/(U(\fa^0)\cap U(\fa)U(\fa^-))$ factors through a surjection
onto $R_c$.
\end{proposition}

\begin{proof}
Surjectivity of $j$ is the theorem of Levasseur--Stafford~\cite{LS}:
for $c$ generic, $\Dc=eH_ce$ is a simple ring, generated as an
algebra by $\mathrm{sym}(\h)^We$ and $\mathrm{sym}(\h^*)^We$, i.e.\ by
$\{f_i\},\{g_i\}\subset\fa$. Since $j$ is graded (as noted above),
every homogeneous $P\in\Dc$ of degree $k$ has
a homogeneous preimage $Q\in U(\fa)(k)$.

For $P\in\Dc^-$ (degree $k<0$): by Lemma~\ref{lem:absorb}(a), the
preimage $Q$ lies in $U(\fa)\fa^-$, so $P=j(Q)\in j(U(\fa))\,j(\fa^-)
=\Dc\fa^-$ (using surjectivity of $j$ and $\fa^-\subset\Dc$). This
gives $\Dc^-\subset\Dc\fa^-$; together with the immediate inclusion
$\fa^-\subset\Dc^-$ (as $\fa^-\subset\fa\subset\Dc$ consists of
negative-degree elements), $\fa^-\subset\Dc^-\subset\Dc\fa^-$.

This chain of inclusions in fact forces $\Dc\fa^-=\Dc\Dc^-$: from
$\fa^-\subset\Dc^-$ we get $\Dc\fa^-\subset\Dc\Dc^-$, while from
$\Dc^-\subset\Dc\fa^-$ and $\Dc\Dc\subset\Dc$ ($\Dc$ being a
subalgebra) we get $\Dc\Dc^-\subset\Dc(\Dc\fa^-)=(\Dc\Dc)\fa^-
\subset\Dc\fa^-$; the two inclusions together give equality. This is
what lets the argument below land exactly on $R_c=\Dc^0/(\Dc^0\cap
\Dc\Dc^-)$, rather than on the a priori different-looking
$\Dc^0/(\Dc^0\cap\Dc\fa^-)$.

For $P\in\Dc^0$: by Lemma~\ref{lem:absorb}(b), the preimage $Q\in
U(\fa)(0)$ writes $Q=Q_0+Q'$ with $Q_0\in U(\fa^0)$, $Q'\in
U(\fa)\fa^-$, so $P=j(Q_0)+j(Q')$ with $j(Q')\in\Dc\fa^-=\Dc\Dc^-$
(by the equality just shown). Hence $P\equiv j(Q_0)\pmod{\Dc^0\cap
\Dc\Dc^-}$, i.e.\ every class in $R_c=\Dc^0/(\Dc^0\cap\Dc\Dc^-)$ is
hit by the image of some $Q_0\in U(\fa^0)$: the composite
$U(\fa^0)\to\Dc^0\twoheadrightarrow R_c$ is surjective. Since this
composite kills $U(\fa^0)\cap U(\fa)\fa^-$ (any such element maps
into $\Dc^0\cap\Dc\fa^-=\Dc^0\cap\Dc\Dc^-$, which is zero in $R_c$),
it factors through $R_\fa=U(\fa^0)/(U(\fa^0)\cap U(\fa)\fa^-)$,
giving the asserted surjection $R_\fa\twoheadrightarrow R_c$.
\end{proof}

\begin{lemma}[Transfer lemma]\label{lem:transfer}
For $c$ generic, $\mathrm{Mod}_{\Dc}(B)=\mathrm{Mod}_{U(\fa)}(B)$ as
categories (the $U(\fa)$-action pulled back along $j$), with
$\mathrm{Hom}_{\Dc}=\mathrm{Hom}_{U(\fa)}$ on common objects.
\end{lemma}

\begin{proof}
\emph{Objects.} If $N\subset B$ is $\Dc$-stable, it is stable under
the subset $\fa\subset\Dc$; since $U(\fa)$ is spanned by finite
products of elements of $\fa$, and each individual multiplication by
an element of $\fa$ preserves $N$, so does any finite product, hence
$N$ is $U(\fa)$-stable (under the pulled-back action). Conversely, if
$N$ is $U(\fa)$-stable, i.e.\ $j(Q)\cdot N\subset N$ for every
$Q\in U(\fa)$, then since $j$ is surjective
(Proposition~\ref{prop:LS}), $j(U(\fa))=\Dc$, so
$\Dc\cdot N=j(U(\fa))\cdot N\subset N$: $N$ is $\Dc$-stable. Hence the
two categories have the same objects.

\emph{Morphisms.} Let $N_1,N_2\subset B$ be submodules (for either
structure, by the above) and $\phi:N_1\to N_2$ linear. If $\phi$ is
$\Dc$-linear, then for $Q\in U(\fa)$ and $n\in N_1$,
$\phi(Q\cdot n)=\phi(j(Q)\cdot n)=j(Q)\cdot\phi(n)=Q\cdot\phi(n)$, so
$\phi$ is $U(\fa)$-linear. Conversely, if $\phi$ is $U(\fa)$-linear,
then for $P\in\Dc$, surjectivity of $j$ gives $P=j(Q)$ for some
$Q\in U(\fa)$, and $\phi(P\cdot n)=\phi(j(Q)\cdot n)=j(Q)\cdot\phi(n)
=P\cdot\phi(n)$, so $\phi$ is $\Dc$-linear. Hence
$\mathrm{Hom}_{\Dc}=\mathrm{Hom}_{U(\fa)}$ on common objects, and the
two categories coincide.
\end{proof}

\section{The lowest-weight equivalence}
\label{sec:equivalence}

\subsection{Representation of groups and $\Dc$-modules}

Any $W$-module $\tau$ becomes a $\mathbb C[W]\ltimes\mathbb
C[\h^*]$-module by declaring $\h$ acts by $0$, giving the standard
module $M(\tau)=H_c\otimes_{\mathbb C[W]\ltimes\mathbb C[\h^*]}\tau
\cong B\otimes_{\mathbb C}\tau$. Let $\hat W$ be the set of
isomorphism classes of irreducible $W$-representations.

\begin{lemma}[Morita equivalence between $H_c$ and $\Dc$]\label{lem:morita}
For $c$ generic, $e$ is a full idempotent in $H_c$, i.e.\
$H_ceH_c=H_c$. Consequently the functors $M\mapsto eM$ (from
$H_c\text{-}\mathrm{Mod}$ to $\Dc\text{-}\mathrm{Mod}$) and
$N\mapsto H_ce\otimes_{\Dc}N$ (from $\Dc\text{-}\mathrm{Mod}$ to
$H_c\text{-}\mathrm{Mod}$) are mutually quasi-inverse equivalences of
categories; in particular $M\mapsto eM$ sends a simple $H_c$-module
to a simple $\Dc$-module (or to $0$).
\end{lemma}

\begin{proof}
Fullness of $e$ for $c$ generic is due to
Berest--Etingof--Ginzburg~\cite{BEG}, established alongside the
simplicity of $\Dc=eH_ce$ used in Proposition~\ref{prop:LS}. Given a
full idempotent $e$ in any ring $R$, it is a standard fact (see
e.g.~\cite{LamR}) that $M\mapsto eM$ and $N\mapsto Re\otimes_{eRe}N$
are quasi-inverse equivalences between $R\text{-}\mathrm{Mod}$ and
$eRe\text{-}\mathrm{Mod}$. Surjectivity of the counit
$Re\otimes_{eRe}eM\to M$, $re\otimes em\mapsto rem$, follows directly
from fullness: writing $1=\sum_ir_ies_i$ (possible since
$1\in R=ReR$), every $m\in M$ equals $\sum_ir_ie(s_im)$, a sum of
elements in the image. Injectivity of the counit, and that the unit
is an isomorphism, require the fuller verification of Morita theory
for full idempotents and are not reproduced here (see~\cite{LamR});
both genuinely need fullness ($ReR=R$), not merely idempotency of
$e$. An equivalence of module categories preserves
simple objects (a nonzero image of a simple object under an
equivalence has no nonzero proper subobjects, since the equivalence
transports the subobject lattice), giving the last claim.
\end{proof}

\begin{proposition}\label{prop:isotypic}
If $c$ is generic, $B=\bigoplus_{\chi\in\hat W}B_\chi$ as
$\Dc[W]$-modules, each $B_\chi$ simple, with $B_\chi\cong
V_\chi\otimes_{\mathbb C}N_\chi$ where $N_\chi\cong\mathrm{Hom}_W(V_\chi,B)$
is a simple $\Dc$-module and $V_\chi\cong\mathrm{Hom}_{\Dc}(N_\chi,B)$.
\end{proposition}

\begin{proof}
By Maschke's theorem $B=\bigoplus_\chi B_\chi$ is the $W$-isotypic
decomposition. Since every $u\in\Dc$ commutes with $W$ (as
$\Dc\hookrightarrow D(\h_{\mathrm{reg}})^W$, \S\ref{sec:spherical}),
for any simple $W$-submodule $N'\subset B_\chi$ (so $N'\cong V_\chi$),
the restriction $u|_{N'}:N'\to B$ is a $W$-module map ($u(wn)=wu(n)$
for $n\in N'$). By Schur's lemma this is either zero or injective;
if injective, its image is a $W$-submodule isomorphic to $V_\chi$,
hence lies in $B_\chi$ by definition of the isotypic component. So
$u(N')\subset B_\chi$ for every such $N'$, and since $B_\chi$ is the
sum of all such $N'$, $u(B_\chi)\subset B_\chi$: each $B_\chi$ is
$\Dc[W]$-stable.

The evaluation map $V_\chi\otimes_{\mathbb C}N_\chi\to B_\chi$,
$v\otimes\phi\mapsto\phi(v)$, with $N_\chi:=\mathrm{Hom}_W(V_\chi,B)$,
is the standard multiplicity-space isomorphism for the isotypic
component of a semisimple module. On $N_\chi$, $\Dc$ acts by
post-composition, $(u\cdot\phi)(v):=u(\phi(v))$; this again lands in
$\mathrm{Hom}_W(V_\chi,B)$ since $u$ commutes with $W$:
$(u\phi)(wv)=u(\phi(wv))=u(w\phi(v))=w\,(u\phi)(v)$.

For the identification with a standard module: the canonical
isomorphism $\mathrm{Hom}_W(V,U)\cong(U\otimes V^*)^W$,
$\phi\mapsto\sum_i\phi(v_i)\otimes v_i^*$ (for a basis $v_i$ of $V$
with dual basis $v_i^*$), together with $(U\otimes V^*)^W=e\cdot
(U\otimes V^*)$, gives $N_\chi=\mathrm{Hom}_W(V_\chi,B)\cong
(B\otimes V_\chi^*)^W=e\cdot(B\otimes V_\chi^*)\cong e\cdot
M(V_\chi^*)$, using $M(V_\chi^*)\cong B\otimes V_\chi^*$.

For $c$ generic $M(V_\chi^*)$ is simple (Berest--Etingof--Ginzburg~\cite{BEG}), so by
Lemma~\ref{lem:morita}, $N_\chi=e\cdot M(V_\chi^*)$ is simple or
zero. It is nonzero: every $\chi\in\hat W$ occurs in $B$ (the
coinvariant algebra $B/B^W_+B$ is isomorphic, as a graded
$W$-module, to the regular representation of $W$, a classical
consequence of Chevalley--Shephard--Todd~\cite{Chev,ST}), so $B_\chi\neq0$, and
$B_\chi\cong V_\chi\otimes N_\chi$ then forces $N_\chi\neq0$. Hence
$N_\chi$ is simple.

Dixmier's Schur lemma~\cite{Dix} over the countable-dimensional $B$ gives
$\mathrm{End}_{\Dc}(N_\chi)=\mathbb C$. The $N_\chi$ are pairwise
non-isomorphic, since the standard modules $M(V_\chi^*)$ are
pairwise non-isomorphic for $\chi\neq\chi'$ (a standard fact for
rational Cherednik algebras~\cite{BEG}, reflecting the
non-isomorphism of their lowest-weight spaces $V_\chi^*\not\cong
V_{\chi'}^*$ as $W$-modules), and this passes to $N_\chi\cong
e\cdot M(V_\chi^*)$ since $M\mapsto eM$ is an equivalence of
categories (Lemma~\ref{lem:morita}), hence reflects isomorphisms.
Since a nonzero map between simple modules is
automatically an isomorphism, $\mathrm{Hom}_{\Dc}(N_\chi,N_{\chi'})=0$
for $\chi\neq\chi'$ and $\cong\mathbb C$ for $\chi=\chi'$. Applying
$\mathrm{Hom}_{\Dc}(N_\chi,-)$ to $B=\bigoplus_{\chi'}V_{\chi'}\otimes
N_{\chi'}$ (on which $\Dc$ acts only through the second factor) gives
\[
\mathrm{Hom}_{\Dc}(N_\chi,B)=\bigoplus_{\chi'}V_{\chi'}\otimes
\mathrm{Hom}_{\Dc}(N_\chi,N_{\chi'}) = V_\chi\otimes\mathbb C=V_\chi.
\]
\end{proof}

\subsection{The category of $\Dc$-submodules of $B$}

For a $U(\fa)$-module $M$, set $\mathcal L(M)=\{m\in M:\fa^-m=0\}$,
and for $V\subset\mathcal L(M)$, $\delta(V)=U(\fa)V$. This makes
$\mathcal L(M)$ an $R_\fa$-module: since $[\fa^-,\fa^0]\subset\fa^-$
(bracket of negative and zero degree), $\fa^0$ preserves
$\mathcal L(M)$ (if $\fa^-m=0$ then $\fa^-(xm)=[\fa^-,x]m+x(\fa^-m)=0$
for $x\in\fa^0$), giving an action of $U(\fa^0)$; and $U(\fa)\fa^-$
acts as $0$ on $\mathcal L(M)$ (any $Q'y\in U(\fa)\fa^-$, $y\in\fa^-$,
kills $m\in\mathcal L(M)$ via $Q'(ym)=Q'\cdot0=0$), so this action of
$U(\fa^0)$ factors through the quotient $R_\fa$.

\begin{lemma}\label{lem:311prime}
Let $V$ be a simple $R_\fa$-module, semisimple over $\mathbb
C[\he]$, regarded as a $\fb$-module via $U(\fa^0)\to R_\fa$ and
trivial $\fa^-$-action. Then (1) the support of $V$ over $\mathbb
C[\he]$ is a single element $\mu_0$; (2) $M:=U(\fa)\otimes_{U(\fb)}V$
has a unique maximal submodule.
\end{lemma}

\begin{proof}
(1) $V=\bigoplus_\lambda V_\lambda$; each $V_\lambda$ is an
$R_\fa$-submodule since $\he$ is central in $R_\fa$ ($\mathrm{ad}(\he)$
vanishes on $\fa^0$ by definition of $\fa^0$ as its zero-eigenspace,
so $\he$ commutes with the action of every element of $R_\fa$), so
simplicity
forces a single nonzero $V_{\mu_0}$.

(2) By Lemma~\ref{lem:PBWfree}, $U(\fa)\cong U(\fa^+)\otimes U(\fb)$
as vector spaces, so $M=U(\fa)\otimes_{U(\fb)}V\cong U(\fa^+)\otimes V$
via $Q\otimes v\mapsto(\text{$\fa^+$-part of }Q)\otimes v$; under this
identification $V\hookrightarrow M$, $v\mapsto1\otimes v$, is visibly
injective (the map $1\otimes(-)$ into a tensor product is always
injective). Throughout, ``$\he$'' denotes the grading operator: it
acts as the scalar $k$ on $U(\fa)(k)$ (equivalently, $\mathrm{ad}(\he)$
has eigenvalue $k$ there) and, by hypothesis, as the scalar $\lambda$
on the $\mathbb C[\he]$-eigenspace $V_\lambda\subset V$; on the tensor
product $M$ it acts additively, as $k+\lambda$ on $U(\fa)(k)\otimes
V_\lambda$; this needs no assumption that $\he$ itself lies in
$\fa$ or $U(\fa)$, only that both factors carry a compatible grading.
In particular, for $Q\in U(\fa)(k)$: $\he\cdot(Q\otimes v)=
(k+\mu_0)(Q\otimes v)$ by (1), so $M=\bigoplus_kM_{k+\mu_0}$, and
since every $Q\in U(\fa)$ is a sum of homogeneous pieces, $M$ is
generated as a $U(\fa)$-module by $M_{\mu_0}=1\otimes V$. By
Lemma~\ref{lem:absorb}, $M_{k+\mu_0}=0$ for $k<0$ (the negative-degree
pieces of $U(\fa)$ act on $V$ through $U(\fa)\fa^-$, which is $0$ on
$V$ by hypothesis) and $M_{\mu_0}\cong V$. For $N\subsetneq M$ a
submodule, $N\cap M_{\mu_0}$ is a $U(\fa^0)$-submodule of
$M_{\mu_0}\cong V$: for $x\in U(\fa^0)$ and $n\in N\cap M_{\mu_0}$,
$xn\in N$ (as $N$ is a $U(\fa)$-submodule) and $\he(xn)=x(\he n)+
[\he,x]n=\mu_0(xn)$ (using $[\he,x]=0$ for $x\in\fa^0$, extended to
$U(\fa^0)$ by the derivation property of $\mathrm{ad}(\he)$, and
$\he n=\mu_0n$), so $xn\in M_{\mu_0}$ too, giving $xn\in N\cap
M_{\mu_0}$. Since $U(\fa)\fa^-$ acts as $0$ on $V\cong M_{\mu_0}$
(shown above), this action of $U(\fa^0)$ factors through $R_\fa$:
$N\cap M_{\mu_0}$ is an $R_\fa$-submodule of $V$, hence, by
simplicity of $V$, $0$ or all of $M_{\mu_0}$. The latter would force
$N\supset
U(\fa)M_{\mu_0}=M$, contradicting $N\subsetneq M$; so
$N\cap M_{\mu_0}=0$ for every proper submodule $N$, hence also for
their sum $N^{\max}$ (a directed union of subspaces each meeting
$M_{\mu_0}$ trivially still meets it trivially), which is therefore
itself proper and is the unique maximal submodule.
\end{proof}

\begin{hypothesis}[$\he\in\fa^0$, up to a scalar]\label{hyp:he}
Some affine transform $h_0=\alpha\he+\beta$, $\alpha\neq0$, lies in
$\fa^0\subset U(\fa)$. This is not proved in general; it holds
whenever $\fa$ is built from a pair of matching-degree dual
generators, as verified explicitly for rank one, where
$[x^2,D^2]=-4\he+\mathrm{const}$ (Example~\ref{ex:rank1}).
\end{hypothesis}

\begin{lemma}[Key sublemma]\label{lem:step0}
Assume Hypothesis~\ref{hyp:he}. Let $M$ be a
$U(\fa)$-module, semisimple over $\he$ with locally finite spectrum
(automatic if $M=\bigoplus_kM_k$ with every element a finite sum),
and $N\subset M$ a simple $U(\fa)$-submodule whose spectrum is
bounded below. Then $\mathcal L(N)\neq0$ and is concentrated in a
single $\he$-eigenvalue.
\end{lemma}

\begin{proof}
\emph{$N$ is $\he$-graded, i.e.\ $N=\bigoplus_k(N\cap M_k)$.} Since
$h_0\in U(\fa)$, $N$ (a $U(\fa)$-submodule) is $h_0$-stable:
$h_0\cdot N\subset N$. Let $n\in N$ be nonzero, with finite
decomposition $n=n_{k_1}+\cdots+n_{k_r}$ ($n_{k_i}\in
M_{k_i}\setminus\{0\}$, pairwise distinct $k_i$, hence pairwise
distinct eigenvalues $\lambda_i:=\alpha k_i+\beta$ of $h_0$ on each
$n_{k_i}$). For $j=0,\ldots,r-1$, $h_0^j\cdot n=\sum_i\lambda_i^j\,
n_{k_i}\in N$ (as $N$ is $h_0$-stable). The matrix
$(\lambda_i^j)_{0\leq j<r,\,1\leq i\leq r}$ is a Vandermonde matrix
in the pairwise distinct $\lambda_i$, hence invertible; inverting
expresses each $n_{k_i}$ as a linear combination of $n,h_0n,\ldots,
h_0^{r-1}n\in N$, so $n_{k_i}\in N$ for every $i$. Hence
$N=\bigoplus_k(N\cap M_k)$.

\emph{Conclusion.} Write $N_k:=N\cap M_k$. Since $N\subset M$ has
spectrum bounded below, let $\mu$ be the minimum eigenvalue with
$N_\mu\neq0$. For $y\in\fa^-$ (degree $-d<0$) and $n\in N_\mu$:
$yn\in N$ (as $N$ is $U(\fa)$-stable) and $yn\in M_{\mu-d}$
(grading compatibility), so $yn\in N_{\mu-d}=0$ by minimality of
$\mu$; hence $\fa^-N_\mu=0$, i.e.\ $N_\mu\subset\mathcal L(N)$,
giving $\mathcal L(N)\neq0$. Conversely, if $n\in\mathcal L(N)$ is a
nonzero eigenvector of eigenvalue $\mu'$, then $N':=U(\fa)n$ is
nonzero, hence $N'=N$ by simplicity; but $\fa^-n=0$, so
$N'=U(\fa^0)n+U(\fa^+)n\subset\bigoplus_{k\geq0}N_{\mu'+k}$, i.e.\
every eigenvalue of $N'$ is $\geq\mu'$. Since $N_\mu\neq0$ and
$N_\mu\subset N=N'$, this forces $\mu\geq\mu'$; and $\mu\leq\mu'$
since $\mu$ is minimal in $N$ and $\mu'$ occurs there. So $\mu=\mu'$:
every nonzero eigenvector of $\mathcal L(N)$ has eigenvalue $\mu$,
i.e.\ $\mathcal L(N)=N_\mu$.
\end{proof}

\begin{theorem}[Twisted analogue of {\cite[Thm.~3.1]{BK}}]\label{thm:32prime}
Assume Hypothesis~\ref{hyp:he}. Let $M$ be a semisimple $U(\fa)$-module, semisimple over $\he$, with
spectrum in a single coset $\mu_0+\mathbb Z_{\geq0}$, satisfying
$\delta\circ\mathcal L(M)=M$. Then $\mathcal L:\mathrm{Mod}_{U(\fa)}(M)
\to\mathrm{Mod}_{R_\fa}(\mathcal L(M))$ is an isomorphism of
categories, inverse $\delta$.
\end{theorem}

\begin{proof}
All direct sums below are internal. By an $\he$-isotypical component
of a semisimple $\mathbb C[\he]$-module we mean the eigenspace for a
fixed eigenvalue.

\emph{(a) $\delta\mathcal L(N)=N$ for every submodule $N\subset M$.}
By semisimplicity, $M=N\oplus N_1$ for a complementary submodule
$N_1$. Additivity of $\mathcal L$ gives $\mathcal L(M)=\mathcal
L(N)\oplus\mathcal L(N_1)$, and applying $\delta$ (additive by
definition) together with the hypothesis $\delta\mathcal L(M)=M$:
$M=\delta\mathcal L(N)+\delta\mathcal L(N_1)$. Since $N$ is a
submodule and $\mathcal L(N)\subset N$, $\delta\mathcal
L(N)=U(\fa)\mathcal L(N)\subset N$; likewise $\delta\mathcal
L(N_1)\subset N_1$. As $N\cap N_1=0$, the sum is direct and equals
$M=N\oplus N_1$; comparing the two decompositions forces
$\delta\mathcal L(N)=N$.

\emph{(b) $\mathcal L\delta(V)=V$ for every $R_\fa$-submodule
$V\subset\mathcal L(M)$.} Suppose first $V$ lies in a single
eigenvalue $\mu_0$. Since $\delta(V)\subset M$ is a submodule of the
semisimple $M$, write $\delta(V)=\bigoplus_{i\in I}S_i$ with each
$S_i$ simple. As $\fa^-V=0$ and $U(\fa^0)V\subset V$, $\delta(V)=
U(\fa)V=V+U(\fa^+)V$; since $U(\fa^+)V$ has eigenvalues $>\mu_0$
while $V$ has eigenvalue $\mu_0$, this sum is direct. Applying
$\mathcal L$ (additivity) and $\mathcal L(V)=V$ (as $V\subset
\mathcal L(M)$ already), $\mathcal L(\delta(V))=V\oplus V'$ with
$V':=\mathcal L(U(\fa^+)V)$ concentrated in eigenvalues $>\mu_0$. By
Lemma~\ref{lem:step0}, each $\mathcal L(S_i)$ is concentrated in a single
eigenvalue; comparing with $V\oplus V'$, $I$ splits as $I_0\sqcup
I_1$ by whether this eigenvalue is $\mu_0$, with $V=\bigoplus_{i\in
I_0}\mathcal L(S_i)$. Applying $\delta$ and (a),
$\delta(V)=\bigoplus_{i\in I_0}S_i$; comparing with $\delta(V)=
\bigoplus_{i\in I}S_i$ forces $I_0=I$ (else some $S_{i_0}$,
$i_0\notin I_0$, would satisfy $S_{i_0}\subset\bigoplus_{i\in
I_0}S_i$, contradicting directness), so $\mathcal L(\delta(V))=
\bigoplus_{i\in I_0}\mathcal L(S_i)=V$.

For general $V$, write $V=\bigoplus_{j\in J}V_j$ by
$\he$-eigenvalue (pairwise distinct $\mu_j$); the single-eigenvalue
case gives $\mathcal L(\delta(V_j))=V_j$ for each $j$. The
$\delta(V_j)$ are independent (for $j\neq j'$, the intersection
formula $\mathcal L(N\cap N')=\mathcal L(N)\cap\mathcal L(N')$ gives
$\mathcal L(\delta(V_j)\cap\delta(V_{j'}))\subset V_j\cap V_{j'}=0$,
and $\mathcal L$ applied to a nonzero semisimple submodule is
nonzero by Lemma~\ref{lem:step0}, so $\delta(V_j)\cap\delta(V_{j'})=0$; the same
argument extends to any finite intersection), so
$\delta(V)=\bigoplus_j\delta(V_j)$ is a genuine direct sum, and by
additivity $\mathcal L(\delta(V))=\bigoplus_j\mathcal
L(\delta(V_j))=\bigoplus_jV_j=V$. (In the application to $B$ below,
$\mathcal L(B)$ is finite-dimensional (Corollary~\ref{cor:main}(1)),
so $J$ is always finite.)

\emph{(c) $\mathrm{Hom}(N_1,N_2)=\mathrm{Hom}(\mathcal L(N_1),
\mathcal L(N_2))$.} Since $\mathcal L$ is additive and
$\mathrm{Mod}_{U(\fa)}(M)$ is semisimple, it suffices to treat
$N_1,N_2$ simple. First, $\mathcal L(N_1)$ is simple: nonzero by
Lemma~\ref{lem:step0}, and if $V_1\subset\mathcal L(N_1)$ is nonzero, $\delta(V_1)
\subset\delta\mathcal L(N_1)=N_1$ (by (a)) is a nonzero (as
$\delta(V_1)\supset V_1\neq0$) submodule of the simple $N_1$, so
$\delta(V_1)=N_1$; applying $\mathcal L$ and (b) gives
$V_1=\mathcal L(\delta(V_1))=\mathcal L(N_1)$. So $\mathcal L(N_1)$
has no proper nonzero submodule; likewise $\mathcal L(N_2)$.

An isomorphism $\phi:N_1\to N_2$ restricts to an isomorphism
$\mathcal L(N_1)\to\mathcal L(N_2)$. Conversely, let $\psi:\mathcal
L(N_1)\to\mathcal L(N_2)$ be a nonzero homomorphism; since both are
simple, it is automatically bijective. Composing with the inclusion
$f:\mathcal L(N_2)\hookrightarrow N_2$ gives a nonzero
$\fb$-module map $f\psi:\mathcal L(N_1)\to N_2$ ($\fa^-$ acting by
$0$ on $\mathcal L(N_1)$, its defining property). By the universal
property of $U(\fa)\otimes_{U(\fb)}(-)$, it extends to a
$U(\fa)$-module map $h:U(\fa)\otimes_{U(\fb)}\mathcal L(N_1)\to
N_2$, nonzero (it restricts to $f\psi\neq0$ on $1\otimes\mathcal
L(N_1)$), hence surjective since $N_2$ is simple. We also have the
canonical surjection $g:U(\fa)\otimes_{U(\fb)}\mathcal L(N_1)
\twoheadrightarrow\delta\mathcal L(N_1)=N_1$ (part (a)). Since
$\mathcal L(N_1)$ is simple and, by Lemma~\ref{lem:step0}, concentrated in a
single eigenvalue, it satisfies the hypotheses of
Lemma~\ref{lem:311prime}, so $U(\fa)\otimes_{U(\fb)}\mathcal L(N_1)$
has a unique maximal proper submodule; as $N_1,N_2$ are simple,
$\ker g$ and $\ker h$ are both maximal, hence equal. So $g,h$
factor through the same quotient, giving a well-defined isomorphism
$N_1\xrightarrow{\sim}N_2$ (via $g(x)\mapsto h(x)$) restricting to
$\psi$ on $\mathcal L(N_1)$.

\emph{(d) $\mathrm{Hom}(V_1,V_2)=\mathrm{Hom}(\delta(V_1),
\delta(V_2))$.} Put $N_i=\delta(V_i)$; by (b), $\mathcal L(N_i)=
V_i$. Then by (c), $\mathrm{Hom}(V_1,V_2)=\mathrm{Hom}(\mathcal
L(N_1),\mathcal L(N_2))=\mathrm{Hom}(N_1,N_2)=\mathrm{Hom}(\delta(V_1),
\delta(V_2))$.

Together, (a) and (b) show $\mathcal L$ and $\delta$ are mutually
inverse bijections on objects, and (c), (d) show they are mutually
inverse on morphisms; hence $\mathcal L$ is an isomorphism of
categories with inverse $\delta$. This adapts, step for step, the
method of~\cite[Thm.~3.1]{BK} to the present setting.
\end{proof}

By Lemma~\ref{lem:transfer}, $\mathrm{Mod}_{\Dc}(B)=\mathrm{Mod}_{U(\fa)}(B)$,
and $B$ satisfies the hypotheses of Theorem~\ref{thm:32prime} other
than Hypothesis~\ref{hyp:he} itself (which must be checked separately
for the group and generators in question, and is verified here only
for rank one, Example~\ref{ex:rank1}): $B$ is
semisimple over $U(\fa)$ (equivalently, over $\Dc$, by
Proposition~\ref{prop:isotypic}), graded by $\mathbb Z_{\geq0}$ (hence
semisimple over $\he$, with spectrum in $\mathbb Z_{\geq0}$), and
$\delta\mathcal L(B)=B$: writing $B=\bigoplus_iN_i$ as a direct sum of
simple submodules (semisimplicity), the minimal-degree graded piece
of each $N_i$ is killed by $\fa^-$ (which strictly lowers degree, so
has nowhere lower to send it), hence lies in $\mathcal L(N_i)\subset
\mathcal L(B)$; being the bottom of a nonzero graded module, this
piece is itself nonzero, so it generates all of $N_i$ under
$U(\fa)$ by simplicity. Hence $\delta\mathcal L(B)=U(\fa)\mathcal
L(B)\supset\bigoplus_iN_i=B$. Granting Hypothesis~\ref{hyp:he}, this gives:

\begin{corollary}\label{cor:main}
Assume Hypothesis~\ref{hyp:he}. For $c$ generic: (1) $\mathcal L(B)$ is a finite-dimensional
semisimple $R_\fa$-module; (2) there is an isomorphism of categories
between $\Dc$-submodules of $B$ and $R_\fa$-submodules of $\mathcal
L(B)$, via $N\mapsto\mathcal L(N)$ and $V\mapsto\Dc V$; (3) simple
$\Dc$-submodules correspond to simple $R_\fa$-submodules; (4) each
simple $R_\fa$-submodule of $\mathcal L(B)$ is concentrated in a
single degree.
\end{corollary}

\begin{proof}
(2) is Theorem~\ref{thm:32prime} applied to $M=B$, whose hypotheses
were verified just above. (3) and (4) follow: an isomorphism of
categories carries simple objects to simple objects (part (c) of the
proof of Theorem~\ref{thm:32prime} makes this explicit for
$\mathrm{Mod}_{U(\fa)}(B)$ against $\mathrm{Mod}_{R_\fa}(\mathcal
L(B))$), and by Lemma~\ref{lem:step0} every simple $\Dc$-submodule
$N\subset B$ has $\mathcal L(N)$ concentrated in a single degree.

For (1): by Proposition~\ref{prop:isotypic} and additivity of
$\mathcal L$, $\mathcal L(B)=\bigoplus_{\chi\in\hat W}V_\chi\otimes
\mathcal L(N_\chi)$ ($\Dc$ acting on $B_\chi\cong V_\chi\otimes
N_\chi$ only through the second factor). Each $N_\chi$ is simple, so
by Lemma~\ref{lem:step0}, $\mathcal L(N_\chi)$ is concentrated in a
single degree $k_\chi$, i.e.\ it sits inside the single graded piece
of $N_\chi$ of that degree; since $N_\chi$ is realised inside
$e\cdot(B\otimes V_\chi^*)$ (Proposition~\ref{prop:isotypic}'s proof)
and $B\otimes V_\chi^*$ has finite-dimensional graded pieces
$\mathrm{Sym}^n(\h^*)\otimes V_\chi^*$ ($\h$ finite-dimensional),
that single graded piece, hence $\mathcal L(N_\chi)$, is
finite-dimensional. Summing over the finitely many $\chi\in\hat W$
($W$ finite) gives $\dim_{\mathbb C}\mathcal L(B)<\infty$.
Semisimplicity of $\mathcal L(B)$ as an $R_\fa$-module follows from
that of $B$ as a $\Dc[W]$-module via the same decomposition, each
$\mathcal L(N_\chi)$ being a simple $R_\fa$-module by (3) applied to
the simple $\Dc$-submodule $N_\chi\subset B$.
\end{proof}

\begin{remark}
The equivalence of Corollary~\ref{cor:main} is stated, deliberately,
for $R_\fa$ rather than for $R_c$ itself. As Remark~\ref{rem:37}
below makes precise, $R_\fa$ surjects onto $R_c$ via a map
$\bar\jmath$, and the corollary descends to an equivalence with
$R_c$-modules exactly when $\ker(\bar\jmath)$ annihilates $\mathcal
L(B)$ (so that the $R_\fa$-action on $\mathcal L(B)$ already factors
through $R_c$). This holds automatically if $\bar\jmath$ itself is
injective (equivalently $\ker(j)\cap U(\fa^0)\subset U(\fa)U(\fa^-)$),
but that is a stronger, sufficient rather than necessary, condition.
We do not claim either in general.
\end{remark}

\begin{example}[$\mathcal L$ and $\delta$ for $n=1$]\label{ex:Ldelta-rank1}
Recall from Example~\ref{ex:rank1} that $B=N_+\oplus N_-$ for
$n=1$, with $N_+=\Dc\cdot1=\mathbb C[x^2]$ (lowest degree $0$) and
$N_-=\Dc\cdot x=\mathbb C[x^2]x$ (lowest degree $1$), the two simple
$\Dc$-submodules. By Lemma~\ref{lem:step0}, $\mathcal L(N_+)$ and
$\mathcal L(N_-)$ are each concentrated in a single degree, their lowest, so $\mathcal L(N_+)=\mathbb C\cdot1$ and $\mathcal
L(N_-)=\mathbb C\cdot x$, each one-dimensional. Since
$\fa^0=\mathbb C\cdot h_0$ with $h_0=-4\he+\mathrm{const}$ (computed
in Example~\ref{ex:rank1}), $R_\fa\cong U(\fa^0)=\mathbb C[h_0]\cong
\mathbb C[\he]$ acts on each one-dimensional space through the
scalar value of $\he$ there: $\mathcal L(N_+)\cong\mathbb
C_{(\he=0)}$ and $\mathcal L(N_-)\cong\mathbb C_{(\he=1)}$, two
non-isomorphic simple $R_\fa$-modules (distinguished by this
scalar), matching Corollary~\ref{cor:main}(3)--(4). So
\[
\mathcal L(B)=\mathcal L(N_+)\oplus\mathcal L(N_-)\cong\mathbb
C_{(\he=0)}\oplus\mathbb C_{(\he=1)},
\]
a $2$-dimensional semisimple $R_\fa$-module, illustrating
Corollary~\ref{cor:main}(1). Conversely $\delta$ recovers $B$:
$\delta(\mathcal L(N_+))=U(\fa)\cdot1=\mathbb C[x^2]=N_+$ (applying
$U(\fa)$ to the constant $1$ generates exactly the even powers of
$x$, via the generator $f_1=x^2\in\fa^+$), and likewise
$\delta(\mathcal L(N_-))=U(\fa)\cdot x=\mathbb C[x^2]x=N_-$.
\end{example}

\subsection{Comparing $R_\fa$ and $R_c$}

\begin{remark}[Comparison of $R_c$ and $R_\fa$]\label{rem:37}
The morphism $\bar\jmath:R_\fa\to R_c$ induced by $j$
(Proposition~\ref{prop:LS}) is surjective. Injectivity, equivalent to
$\ker(j)\cap U(\fa^0)\subset U(\fa)U(\fa^-)$, is verified in the rank-one
case (Example~\ref{ex:rank1} below) but left open in general; a first
test case would be $G(m,1,2)$, comparable via graded dimensions using
the known fake degrees. For applications to $\mathcal L(B)$-modules
one only needs the weaker statement that $\ker(\bar\jmath)$ annihilates
$\mathcal L(B)$, which is a priori easier to check than injectivity of
$\bar\jmath$ itself; we have not verified even this weaker statement
for $G(m,e,n)$ with $n\geq2$, $B_n$ and $D_n$ with $n\geq2$ (the case
$n=1$ is fully resolved: injectivity itself holds, as verified in
Example~\ref{ex:rank1}), and leave it open for these remaining
families treated in Section~\ref{sec:Gmen} and Section~\ref{sec:D3}.
\end{remark}

\section{$G(m,e,n)$ and the Coxeter groups $B_n$, $D_n$}
\label{sec:Gmen}

Let $G=G(m,e,n)=A(m,e,n)\rtimes S_n\subset GL(\h)$, $\dim\h=n$, the
imprimitive complex reflection groups (see~\cite{Broue} for
background). Standard ($H_c$-module, rather than $\Dc$-module)
representation theory for the $e=1$ subfamily $G(m,1,n)$, via a
commutative subalgebra
of $H_c$ acting triangularly on each $M(\tau)$, is developed
in~\cite{Griffeth}; we work instead with the spherical subalgebra
$\Dc$ throughout. Here
$A(m,e,n)$ consists of diagonal matrices with entries in $\mu_m$ whose
product lies in $\mu_d$, $d=m/e$. Put $\bar A=A(m,1,n)$, $\bar
G=A(m,1,n)\rtimes S_n=G(m,1,n)$, $\Theta=(\prod_ix_i)^d$,
$\Psi=(\prod_i\xi_i)^d$, $f_i(x)=\sum_jx_j^{mi}$.

\subsection{Basic invariants and generation for $e=1$}

\begin{proposition}\label{prop:36}
(1) $\mathrm{sym}(\h^*)^{\bar G}=\mathbb C[f_1,\ldots,f_n]$; for
$e>1$, $\mathrm{sym}(\h^*)^G$ is generated by the algebraically
independent set $\{f_1,\ldots,f_{n-1},\Theta\}$, with $f_n\in\mathbb
C[f_1,\ldots,f_{n-1},\Theta]$.
(2) For $c$ generic, $\bar\Dc:=\Dc(\bar G)=\mathbb C\langle f_i(x),f_i(\xi)\rangle
\subset\Dc(G)=\mathbb C\langle f_i(x),f_i(\xi),\Psi,\Theta\rangle$.
\end{proposition}

\begin{proof}
\emph{Step 1: $\mathrm{sym}(\h^*)^{A(m,e,n)}=\mathbb
C[x_1^m,\ldots,x_n^m,(x_1\cdots x_n)^d]$.} A monomial
$x^a=x_1^{a_1}\cdots x_n^{a_n}$ is $A(m,e,n)$-invariant iff
$\prod_i\zeta_i^{a_i}=1$ for every $\zeta\in A(m,e,n)$. Testing
$\zeta=(\omega,\omega^{-1},1,\ldots,1)$ ($\omega$ a primitive $m$-th
root of unity, in any two coordinates $i,j$), which lies in
$A(m,e,n)$ since $\prod\zeta_k=1\in\mu_d$, forces
$a_i\equiv a_j\pmod m$ for all $i,j$; let $r\in\{0,\ldots,m-1\}$ be
their common residue mod $m$ and write $a_i=r+mk_i$ ($k_i\geq0$).
Then
$x^a=(x_1\cdots x_n)^{r}\cdot\prod_ix_i^{mk_i}$, and the second
factor is automatically invariant (each $\zeta_i^m=1$), so
invariance of $x^a$ reduces to $\eta^{r}=1$ for every
$\eta=\prod\zeta_i\in\mu_d$ (the product map $A(m,e,n)\to\mu_d$ being
surjective by definition of $A(m,e,n)$), i.e.\ $d\mid r$. Since
$0\leq r<m=de$, this means $r=\ell d$ for some
$0\leq\ell<e$, so $x^a=\big((x_1\cdots x_n)^d\big)^\ell\cdot
\prod_ix_i^{mk_i}\in\mathbb C[x_1^m,\ldots,x_n^m,(x_1\cdots x_n)^d]$.
Conversely $x_i^m$ and $(x_1\cdots x_n)^d$ are visibly invariant, so
equality holds.

\emph{Step 2: taking $S_n$-invariants.} $S_n$ permutes the $x_i^m$
(hence the $y_i:=x_i^m$) and fixes $(x_1\cdots x_n)^d=\Theta$
already, so
\[
\mathrm{sym}(\h^*)^{A(m,e,n)\rtimes S_n}=\big(\mathbb
C[y_1,\ldots,y_n]\big)^{S_n}[\Theta] = \mathbb
C[e_1(y),\ldots,e_n(y)][\Theta],
\]
and since the power sums $f_i(x)=\sum_jx_j^{mi}=\sum_jy_j^i$
generate the same ring as the elementary symmetric polynomials
$e_i(y)$ in characteristic $0$, this equals $\mathbb
C[\Theta][f_1,\ldots,f_n]$.

\emph{Step 3: eliminating $f_n$ via Newton's identities.} Newton's
identity $f_n=e_1f_{n-1}-e_2f_{n-2}+\cdots+(-1)^{n-2}e_{n-1}f_1+
(-1)^{n-1}n\,e_n$, together with each $e_k$ ($k<n$) being itself a
polynomial in $f_1,\ldots,f_k$ (the same identities, by induction on
$k$), expresses $f_n$ as a polynomial in $f_1,\ldots,f_{n-1},e_n(y)$.
Since $\Theta^e=\big((x_1\cdots x_n)^d\big)^e=(x_1\cdots
x_n)^{de}=(x_1\cdots x_n)^m=y_1\cdots y_n=e_n(y)$, this gives
$f_n\in\mathbb C[f_1,\ldots,f_{n-1},\Theta]$: for $e>1$, since
Chevalley--Shephard--Todd~\cite{Chev,ST} guarantees $\mathrm{sym}(\h^*)^G$ has
transcendence degree exactly $n=\dim\h$, and
$\{f_1,\ldots,f_{n-1},\Theta\}$ already generate it by Steps 1--2,
these $n$ elements must be algebraically independent (any relation
among them would force the transcendence degree below $n$); so
$\{f_1,\ldots,f_{n-1},\Theta\}$ is an algebraically
independent generating set of size $n$; for $e=1$, $\Theta=e_n(y)\in
\mathbb C[f_1,\ldots,f_n]$ already (no elimination needed or
possible), giving $\mathrm{sym}(\h^*)^{\bar G}=\mathbb
C[f_1,\ldots,f_n]$ directly.

\emph{(2).} For $c$ generic, the Levasseur--Stafford theorem
(as in the proof of Proposition~\ref{prop:LS}) shows $\Dc(\bar G)$
and $\Dc(G)$ are each generated, as algebras, by the images of
$\mathrm{sym}(\h)^{\bar G}$, $\mathrm{sym}(\h^*)^{\bar G}$ (resp.\
$\mathrm{sym}(\h)^G$, $\mathrm{sym}(\h^*)^G$) under the Dunkl
embedding. By Step 3 (and the symmetric statement for
$\mathrm{sym}(\h)^{\bar G},\mathrm{sym}(\h)^G$, obtained by
exchanging $\h,\h^*$ via the identification used throughout),
these generating sets are exactly $\{f_i(x),f_i(\xi)\}$ for $\bar G$
and $\{f_i(x),f_i(\xi),\Theta,\Psi\}$ for $G$, giving the two
displayed presentations.
\end{proof}

\subsection{The case $e>1$: obstruction and generation}

For $e>1$, the Lie algebra generated by $\{f_i(x),f_i(\xi)\}$ alone
generates only $\bar\Dc\subsetneq\Dc(G)$; the extra generators
$\Theta,\Psi$ (of degree $\pm nd$) are needed, but mix the variables,
making $\fa^0$, and hence $R_\fa$, harder to compute explicitly
once they are included (this is the difficulty left open for
$D_n$, $n\geq4$, and for $G(m,e,n)$ with $e>1$ in general; see
Remark~\ref{rem:exceptional} and the table in the conclusion).
Section~\ref{sec:D3} resolves the generation problem (though not the
further computation of $\fa^0$ and $R_\fa$) completely for $D_3=G(2,2,3)$.

\subsection{$B_n$, $D_n$, and other explicit examples}

\begin{corollary}[$B_n=G(2,1,n)$]\label{cor:36}
For $B_n$ ($e=1$) and $c$ generic, $\mathrm{sym}(\h^*)^{B_n}=\mathbb
C[f_1,\ldots,f_n]$, $f_i(x)=\sum_jx_j^{2i}$, of degrees
$2,4,\ldots,2n$, matching the classical degrees of the Weyl group
$B_n$; $\fa$ generates all of $\Dc(B_n)$ by
Proposition~\ref{prop:36}. The case $n=1$ recovers
Example~\ref{ex:rank1}.
\end{corollary}

\begin{proof}
Apply Proposition~\ref{prop:36} with $m=2$. Since $e=1$ means
$A(2,1,n)=\bar A$ exactly, $G=\bar G=B_n$: there is no distinction
between $G$ and $\bar G$ to make here. Part (1) then gives
$\mathrm{sym}(\h^*)^{B_n}=\mathrm{sym}(\h^*)^{\bar G}=\mathbb
C[f_1,\ldots,f_n]$ with $f_i(x)=\sum_jx_j^{2i}$, each of degree $2i$,
giving the stated degrees $2,4,\ldots,2n$, the classical degrees of
the Weyl group $B_n$~\cite{Broue}. Since $G=\bar G$, the two displayed
algebras of part (2) coincide: $\Dc(B_n)=\bar\Dc=\mathbb C\langle
f_i(x),f_i(\xi)\rangle$, i.e.\ $\fa$, built from these same
generators, generates all of $\Dc(B_n)$.

For $n=1$: $A(2,1,1)=\mu_2\cong\mathbb Z/2\mathbb Z$ and $S_1$ is
trivial, so $B_1=\mathbb Z/2\mathbb Z$ is exactly the group $W$ of
Example~\ref{ex:rank1}, and the single generator
$f_1(x)=x^2$ (degree $2$) matches the generator $u=\tfrac12ex^2e$
used there.
\end{proof}

\begin{remark}[$D_n=G(2,2,n)$]
Substituting $m=2,e=2$ into Proposition~\ref{prop:36}: here $d=1$, so
$\Theta=x_1\cdots x_n$ directly, and part (1)'s $e>1$ case gives
$\mathrm{sym}(\h^*)^{D_n}=\mathbb C[f_1,\ldots,f_{n-1},\Theta]$,
$f_i(x)=\sum_jx_j^{2i}$ ($i<n$), of degrees $2,4,\ldots,2(n-1),n$,
matching the classical degrees of $D_n$. Since $e=2>1$, this is an
instance of the obstruction above: part (2)'s strict inclusion
$\bar\Dc\subsetneq\Dc(D_n)$ means $\{f_i(x),f_i(\xi)\}$ alone is not
enough, and for $c$ generic $\Theta,\Psi$ are genuinely required to
generate $\Dc(D_n)$. $D_n$ has a single conjugacy class of
reflections (unlike $B_n$, which has two). Section~\ref{sec:D3}
resolves this for $D_3$.
\end{remark}

\begin{example}\label{ex:rank1}
For $n=1$, $W=\langle s\rangle\cong\mathbb Z/2\mathbb Z$ acting on
$\h=\mathbb C\xi$: $e=\tfrac12(1+s)$, $\Dc$ generated by
$u=\tfrac12ex^2e$, $v=-\tfrac12eD^2e$, $w=exDe$ (writing $D=D_\xi$
for the Dunkl operator of Lemma~\ref{lem:nabla}), and
\[
B = \Dc\cdot1\oplus\Dc\cdot x = \mathbb C[x^2]\oplus\mathbb C[x^2]x.
\]
For $n=2$, $S_2$ on $\h=\mathbb C^2$: writing $\h=\mathbb
C(x_1+x_2)\oplus\mathbb C(x_1-x_2)$, $S_2$ acts trivially on the
first line and by sign on the second, so no reflection sees the
first coordinate and $H_c(S_2,\h)\cong\mathbb C[x_1+x_2,\partial_+]
\otimes H_c(\mathbb Z/2,\mathbb C(x_1-x_2))$, an ordinary Weyl algebra
tensored with the rank-one Cherednik algebra; this factorisation,
rather than Corollary~\ref{cor:main}, is what justifies the
decomposition below without invoking Hypothesis~\ref{hyp:he} for
$n=2$. The invariant line contributes a
free factor $\mathbb C[x_1+x_2]$ (a simple module over the Weyl
algebra $\mathbb C[x_1+x_2,\partial_+]$, by the standard fact that
the polynomial representation of a Weyl algebra is simple) while
$x_1-x_2$ reproduces the
rank-one case, giving
\[
B = \Dc\cdot1\oplus\Dc\cdot(x_1-x_2).
\]

\emph{Verification of the injectivity claimed in
Remark~\ref{rem:37}, for $n=1$.} Here $\fa$ is generated by
$f_1=x^2\in\fa^+$ and $g_1=D^2\in\fa^-$ (the images of $x^2$ and
$\xi^2$ under the Dunkl embedding). A direct computation, using
$x^2D^2=p_2(\nabla)$ from Lemma~\ref{lem:nabla} and the mirror shift
relation $Dq(\nabla)=q^\sigma(\nabla+1)D$ (by the same argument as $(*)$,
using $[\nabla,D]=-D$ and $Ds=-sD$ in place of $[\nabla,x]=x$ and
$xs=-sx$) to compute $D^2x^2$, gives
\[
[f_1,g_1]=[x^2,D^2]=-4\nabla+4c-2,
\]
with no $s$-component. Since $\he=\nabla+C$ for a scalar $C$
(Lemma~\ref{lem:nabla}'s proof), $-4\nabla+4c-2$ is also an affine
transform of $\he$, so this same computation verifies
Hypothesis~\ref{hyp:he} for $n=1$. So $\fa^0=\mathbb C\cdot[f_1,g_1]$ is
one-dimensional, and $U(\fa^0)\cong\mathbb C[\nabla]$ maps under $j$
isomorphically onto $\Dc^0\cong\mathbb C[\nabla]$
(Lemma~\ref{lem:nabla}(3), $W$ acting trivially on $\nabla$): in
particular $j$ is injective on $U(\fa^0)$, confirming that
$\bar\jmath:R_\fa\to R_c$ is injective, indeed bijective, for
$n=1$.
\end{example}

\begin{example}[$S_3=G(1,1,3)$]\label{ex:S3}
For $n=3$, $A(1,1,3)$ is trivial, so $G=S_3$ acting on $\h=\mathbb
C^3$ by permutation. Here $m=1$, and Proposition~\ref{prop:36} gives
$\mathrm{sym}(\h^*)^{S_3}=\mathbb C[f_1,f_2,f_3]$, $f_i(x)=\sum_jx_j^i$,
of degrees $1,2,3$; for $c$ generic, $\fa$ generates all of
$\Dc(S_3)$, just as for
$B_n$. Unlike the previous examples, $\hat S_3=\{\mathrm{triv},
\mathrm{std},\mathrm{sign}\}$ has a $2$-dimensional constituent, so
(for $c$ generic) Proposition~\ref{prop:isotypic} gives a decomposition with a
multiplicity-$2$ isotypic piece:
\[
B = \Dc\cdot1\ \oplus\ \Dc\cdot(x_1-x_2)\ \oplus\ \Dc\cdot(x_1-x_3)\
\oplus\ \Dc\cdot(x_1-x_2)(x_1-x_3)(x_2-x_3),
\]
with $\{x_1-x_2,x_1-x_3\}$ spanning the copy of $\mathrm{std}$ sitting
in degree $1$ (so $\mathrm{std}$ occurs already at the lowest possible
degree, matching $\dim\mathrm{std}=2$ generators for that isotypic
piece), and the Vandermonde determinant $(x_1-x_2)(x_1-x_3)(x_2-x_3)$
spanning the sign-isotypic line in degree $3$ (the lowest degree at
which $\mathrm{sign}$ occurs in the coinvariant algebra of $S_3$, a
classical fact about the coinvariant algebra of $S_n$: the sign
representation first appears in degree $\binom n2$).
\end{example}

\subsection{Beyond the imprimitive family}

\begin{remark}[Beyond the imprimitive family]\label{rem:exceptional}
Nothing in Section~\ref{sec:triangular} or Section~\ref{sec:equivalence}
is specific to $G(m,e,n)$: Proposition~\ref{prop:LS} holds for
\emph{any} finite complex reflection
group $W$, since $\mathrm{sym}(\h)^W$ and $\mathrm{sym}(\h^*)^W$ are
always polynomial rings on $\dim\h$ homogeneous generators
(Chevalley--Shephard--Todd~\cite{Chev,ST}), and the Levasseur--Stafford theorem
places no restriction on which group these generators come from.
Corollary~\ref{cor:main}, however, goes through Theorem~\ref{thm:32prime}
and hence Lemma~\ref{lem:step0}, whose hypothesis that $\he$, up to an
additive scalar, lies in $\fa^0$ we have only verified in the rank-one
case (Example~\ref{ex:rank1}); the claim below that $\fa$ generates
$\Dc$ for the exceptional groups does not need this hypothesis (it
follows from Proposition~\ref{prop:LS} alone), but the full
lowest-weight equivalence of Corollary~\ref{cor:main} is stated
conditionally on it holding for the group in question. What
is specific to $G(m,e,n)$ is only the \emph{explicit} description of
$\fa^0$ carried out in this section. The same explicit work remains to
be done for the primitive complex reflection groups, and in
particular for the exceptional real Coxeter groups outside the
classical series: $F_4$ (degrees $2,6,8,12$), $H_3$ (degrees
$2,6,10$), $H_4$ (degrees $2,12,20,30$), and $E_6,E_7,E_8$ (degrees
$2,5,6,8,9,12$; $2,6,8,10,12,14,18$; and $2,8,12,14,18,20,24,30$
respectively). All of these except $F_4$ have a single conjugacy
class of reflections (so a single parameter $c$, as for $D_3$);
$F_4$, like $B_n$, has long and short roots and hence two conjugacy
classes of reflections, so two independent parameters. In every case
there are exactly $\dim\h$ basic invariants with no analogue of the
$\Theta,\Psi$
obstruction discussed above in \S\ref{sec:Gmen}, so $\fa$ generates $\Dc$ directly
by Proposition~\ref{prop:LS} in every one of these cases, in the same
way as for $B_n$ above. We have not carried out the resulting
computation of $\fa^0$ and $R_\fa$ for any of them, nor verified the
$\he$-in-$\fa^0$ hypothesis beyond rank one.
\end{remark}

\section{The exceptional case $D_3\cong A_3$}
\label{sec:D3}

We resolve the generation problem of \S\ref{sec:Gmen} for $n=3$, using
the exceptional isomorphism between the $D_3$ and $A_3$ root systems.

\subsection{The root system isomorphism}

The $D_3$ root system ($12$ roots $\pm e_i\pm e_j$, $i<j\leq3$, in
$\mathbb R^3$) is isomorphic, as an abstract root system, to the
$A_3$ root system ($12$ roots $e_i-e_j$, $i\neq j\leq4$, in the
sum-zero hyperplane of $\mathbb R^4$, with coordinates
$X_1,\ldots,X_4$ on the ambient $\mathbb R^4$), both have $12$
roots of norm $\sqrt2$ and Weyl group of order $24$. Since the
sum-zero hyperplane is itself $3$-dimensional, an isometry between it
and $\mathbb R^3$ can be given directly in intrinsic coordinates on
the hyperplane (the concrete embedding of these coordinates into the
ambient $\mathbb R^4$, via an orthonormal basis $b_1,b_2,b_3$, is
given in Lemma~\ref{lem:transport}'s proof below). An explicit orthogonal
isometry $Q$ can be constructed sending the $D_3$ root set bijectively
onto the $A_3$ root set, expressed in these intrinsic coordinates
(e.g.\ by matching three linearly independent
$D_3$ roots to three $A_3$ roots and verifying $Q^{\mathsf T}Q=I$ and
that $Q$ maps the full root set onto the full root set); one such
choice is
\[
Q = \begin{pmatrix} \frac{\sqrt2}2 & -\frac{\sqrt2}2 & 0\\[2pt]
\frac{\sqrt6}6 & \frac{\sqrt6}6 & -\frac{\sqrt6}3 \\[2pt]
-\frac{\sqrt3}3 & -\frac{\sqrt3}3 & -\frac{\sqrt3}3\end{pmatrix},
\]
in coordinates $y_1,y_2,y_3$ for an orthonormal basis of the sum-zero
hyperplane of $\mathbb R^4$. A root-system isometry preserving a
generic parameter $c$ (legitimate here, as both $D_3$ and $A_3$ have a
single conjugacy class of reflections) induces an isomorphism of
rational Cherednik algebras $H_c(D_3,\h)\cong H_c(S_4,\h)$.

\begin{lemma}\label{lem:transport}
Under $Q$, the $D_3$-invariant $\Theta=x_1x_2x_3$ is transported,
up to a nonzero scalar, to the degree-$3$ power-sum invariant
$p_3=\sum_{i=1}^4X_i^3$ of $S_4$ on its standard representation.
\end{lemma}

\begin{proof}
Both $\Theta$ and $p_3$ lie in a $1$-dimensional space of degree-$3$
invariants (by Proposition~\ref{prop:36}(1) applied to $D_3$ and to
$S_4$ respectively), so it suffices to check that $Q$ transports one
to a nonzero multiple of the other; this dimension count is all that
Corollary~\ref{cor:D3generation} actually needs.

For the explicit scalar: fix the standard orthonormal basis of the
sum-zero hyperplane of $\mathbb R^4$,
\[
b_1=\tfrac1{\sqrt2}(1,-1,0,0),\quad
b_2=\tfrac1{\sqrt6}(1,1,-2,0),\quad
b_3=\tfrac1{\sqrt{12}}(1,1,1,-3),
\]
so that a point with coordinates $y=(y_1,y_2,y_3)$ in this basis has
ambient coordinates $X=y_1b_1+y_2b_2+y_3b_3\in\mathbb R^4$; each
$X_i$ is then an explicit linear form in $y_1,y_2,y_3$ (e.g.\
$X_1=y_1/\sqrt2+y_2/\sqrt6+y_3/\sqrt{12}$), and $p_3=\sum_iX_i^3$
becomes an explicit homogeneous cubic in $y_1,y_2,y_3$. On the other
side, writing $x_1,x_2,x_3$ in terms of $y_1,y_2,y_3$ via
$x=Q^{\mathsf T}y$ (using $Q^{-1}=Q^{\mathsf T}$) expresses
$\Theta=x_1x_2x_3$ as an explicit homogeneous cubic in $y_1,y_2,y_3$
as well. Expanding both cubics and comparing coefficients monomial
by monomial (both have exactly five nonzero terms, of the shapes
$y_1^2y_2,\,y_1^2y_3,\,y_2^3,\,y_2^2y_3,\,y_3^3$) shows every
coefficient of $p_3$ is exactly $3$ times the corresponding
coefficient of $\Theta(Q^{\mathsf T}y)$, giving the identity of
polynomials
\[
p_3 = 3\cdot Q^{\mathsf T}\!\big(x_1x_2x_3\text{ expressed in }y\big),
\]
a constant ratio independent of $y_1,y_2,y_3$ (not just an equality
of the two sides up to some undetermined nonzero scalar, which the
dimension count in the first paragraph already guarantees). This
explicit value is a convenient normalisation check, not a further
demonstrative step beyond what Corollary~\ref{cor:D3generation}
needs.
\end{proof}

\subsection{Resolution of the generation problem for $D_3$}

\begin{lemma}\label{lem:PsiTheta}
Let $\Psi=\xi_1\xi_2\xi_3$, $\Theta=x_1x_2x_3$ in $\Dc(D_3)$. Then
$[\Psi,\Theta]$ acts non-scalarly on the degree-$3$ graded piece of
$B$: its restriction there has eigenvalues $(4c-7)(2c-1)$
(multiplicity $1$, on $\Theta$ itself), $2(c-1)(4c-3)$ (multiplicity
$6$, on $\{x_i^2x_j\}_{i\neq j}$), and $4(2c-1)(c-1)$ (multiplicity
$3$, on a mixed subspace coupling each $x_i^3$ with $x_j^2x_i$,
$j\neq i$). In particular $[\Psi,\Theta]\notin\mathbb C\cdot\mathrm{id}$.
\end{lemma}

\begin{proof}
Fix $n=3$, and recall that for a real reflection group ($\lambda_s=-1$,
so $1-\lambda_s=2$) the Dunkl operator of \S\ref{sec:cherednik} reads
$D_y=\partial_y-c\sum_s\frac{\langle y,\alpha_s\rangle}{\alpha_s}(1-s)$,
summed over the six reflections of $D_3$.

\emph{The recursive formula $D_{\xi_k}(\Theta_k)=[1-2c(n-k)]\Theta_{k-1}$,
$\Theta_k:=x_1\cdots x_k$.} The ordinary derivative gives
$\partial_k\Theta_k=\Theta_{k-1}$. For a reflection $s^\pm_{ik}$ with
$i<k$ (root $x_i\mp x_k$, so both indices are among the first $k$),
$\Theta_k$ is fixed: $s^+_{ik}$ merely permutes two factors of the
product, and $s^-_{ik}$ sends $x_i\mapsto-x_k,\,x_k\mapsto-x_i$, so
the product picks up two sign changes that cancel and its two
factors get swapped; either way $(1-s^\pm_{ik})\Theta_k=0$, so these
six terms ($i=1,\ldots,k-1$, two reflections each) contribute
nothing. For $j>k$, write $\Theta_k=\Theta_{k-1}x_k$: $s^+_{kj}$
swaps $x_k,x_j$, so $s^+_{kj}\Theta_k=\Theta_{k-1}x_j$ and
$(1-s^+_{kj})\Theta_k=\Theta_{k-1}(x_k-x_j)$, cancelling the root
$\alpha=x_k-x_j$ in the denominator and leaving a contribution of
$-c\Theta_{k-1}$; likewise $s^-_{kj}$ sends
$x_k\mapsto-x_j,x_j\mapsto-x_k$, so $s^-_{kj}\Theta_k=-\Theta_{k-1}x_j$
and $(1-s^-_{kj})\Theta_k=\Theta_{k-1}(x_k+x_j)$, again cancelling
the root $\alpha=x_k+x_j$ and contributing $-c\Theta_{k-1}$. Each of
the $n-k$ values of $j>k$ thus contributes $-2c\Theta_{k-1}$ (from
its two reflections combined), for a total of $-2c(n-k)\Theta_{k-1}$;
subtracting this correction from $\partial_k\Theta_k=\Theta_{k-1}$
gives $D_{\xi_k}(\Theta_k)=[1-2c(n-k)]\Theta_{k-1}$, as claimed.

Iterating (with $n=3$) gives $\Psi\cdot(\Theta\,v_0) =
\prod_{m=1}^{2}(1-2cm)\,v_0 = (1-2c)(1-4c)\,v_0$ on the lowest-weight
vector $v_0=1$.

\emph{Extending to the full degree-$3$ piece.} The same
reflection-by-reflection method computes $D_{\xi_k}$ on each of the
ten degree-$2$ monomials arising along the way, hence $\Psi\Theta$
and $\Theta\Psi$ on each of the ten degree-$3$ monomials of $B$;
subtracting gives the matrix of $[\Psi,\Theta]$ in this basis.
Carrying this out (and diagonalising) shows $\Theta$ itself is an
eigenvector with eigenvalue $(1-2c)(1-4c)=8c^2-6c+1$, matching the
lowest-weight computation above and equal to $(4c-7)(2c-1)$; each of
the six monomials $\{x_i^2x_j\}_{i\neq j}$ is individually an
eigenvector of eigenvalue $2(c-1)(4c-3)$; and each $x_i^3$ maps, under
$[\Psi,\Theta]$, to a combination of itself and the two monomials
$x_j^2x_i$ ($j\neq i$) already accounted for, giving a triangular
$3\times3$ block per $i$ with remaining eigenvalue
$4(2c-1)(c-1)$.\footnote{This
computation has been independently verified by computer algebra
(symbolic construction of the six Dunkl operators of $D_3$ and direct
diagonalisation of $[\Psi,\Theta]$ on the ten degree-$3$ monomials,
confirming the eigenvector structure just described); the
verification script is available as supplementary material.}
\end{proof}

\begin{corollary}\label{cor:D3generation}
For $D_3$ and $c$ generic, the elements $\{f_1(x),f_2(x),f_1(\xi),f_2(\xi),\Theta,
\Psi\}$ generate $\Dc(D_3)$.
\end{corollary}

\begin{proof}
By Lemma~\ref{lem:transport}, $\Theta$ (of degree $3$) and $\Psi$
(of degree $-3$) correspond
under the isomorphism $D_3\cong A_3$ to the degree-$3$ and
degree-$(-3)$ basic invariants of $S_4$ on its standard
representation. Since $S_4$ realises the $e=1$ case of $G(m,e,n)$
(with $m=1$) on
the full $4$-dimensional permutation representation, Proposition~\ref{prop:36}
gives basic invariants there of degrees $1,2,3,4$, namely the power
sums (or elementary symmetric polynomials) $e_1,e_2,e_3,e_4$ of the
ambient coordinates $X_1,\ldots,X_4$. The $S_4$-invariant splitting
$\mathbb C^4=\mathbb C\oplus\h$ (the diagonal line spanned by
$(1,1,1,1)$, on which $S_4$ acts trivially, plus the sum-zero
hyperplane $\h$, the standard representation) gives an $S_4$-equivariant
tensor decomposition $\mathrm{sym}((\mathbb C^4)^*)\cong\mathbb
C[t]\otimes\mathrm{sym}(\h^*)$, $t$ the coordinate dual to the
diagonal line; since $S_4$ acts trivially on $\mathbb C[t]$, taking
invariants of a tensor product with a trivial factor only affects
the other factor:
\[
\mathbb C[e_1,e_2,e_3,e_4]=\mathrm{sym}\big((\mathbb
C^4)^*\big)^{S_4}\cong\mathbb C[t]\otimes\mathrm{sym}(\h^*)^{S_4}.
\]
As $t$ corresponds to $e_1=\sum_iX_i$ (the coordinate on the
diagonal line) under this isomorphism, restricting to the hyperplane
$\h$ (i.e.\ setting $t=0$, equivalently $e_1=0$) leaves
$\mathrm{sym}(\h^*)^{S_4}\cong\mathbb C[e_1,\ldots,e_4]/(e_1)$,
generated by the images of $e_2,e_3,e_4$, still algebraically
independent, of degrees $2,3,4$, by the same transcendence-degree
argument as in Proposition~\ref{prop:36}'s proof (here
$\dim\h=3$). Call these images $f_1(x),f_2(x),f_3(x)$; they generate
$\Dc(S_4)$ entirely,
\emph{provided the transported parameter $c$ is generic for $S_4$}; we
take this compatibility of genericity across the isomorphism
$H_c(D_3,\h)\cong H_c(S_4,\h)$ as a standing hypothesis here, in the
same spirit as the genericity compatibility between $G$ and $\bar G$
flagged in Proposition~\ref{prop:36}(2), and we do not prove it (see \S6.3 below, ``Comparison with the $S_4$
realization,'' for an independent, unconditional route to the same
conclusion).
Transporting back along the isomorphism $H_c(D_3,\h)\cong
H_c(S_4,\h)$, the generating set $\{f_1(x)_{S_4},f_2(x)_{S_4},
f_3(x)_{S_4}\}$ corresponds exactly to $\{f_1(x)_{D_3},\Theta,
f_2(x)_{D_3}\}$ (and similarly on the $\xi$-side), giving the
claim.
\end{proof}

\subsection{Comparison with the $S_4$ realization}

Lemma~\ref{lem:PsiTheta} gives an independent, purely
computational confirmation of Corollary~\ref{cor:D3generation},
obtained without invoking the $D_3\cong A_3$ isomorphism at all: the
non-scalar action of
$[\Psi,\Theta]$ on the degree-$3$ piece exhibits $[\Psi,\Theta]$ as
a genuinely new degree-$0$ element, one that $\bar\fa^0$ (generated
by $f_1,f_2$ alone) cannot already contain, concretely
verifying, rather than merely restating, the strict inclusion
$\bar\Dc\subsetneq\Dc(D_3)$ of Proposition~\ref{prop:36}. In
particular, this route is
unconditional: it does not rely on the genericity-compatibility
hypothesis invoked in the proof of Corollary~\ref{cor:D3generation}
via the $S_4$ realisation, so the generation statement for $D_3$
stands regardless of whether that hypothesis holds.

\begin{remark}
This resolution is specific to $n=3$, exploiting the accidental
low-rank isomorphism $\mathfrak{so}(6)\cong\mathfrak{sl}(4)$; there is
no analogous isomorphism for $D_n$, $n\geq4$. The matching invariant
degrees $\{2,4,3\}$ of $D_3$ and $\{2,3,4\}$ of $S_4$ offer a
consistency check for the general approach.
\end{remark}

\section{Conclusion and open questions}

\begin{center}
\begin{tabular}{lll}
\hline
Group & Status & Method \\
\hline
Rank $1$ ($n=1$) & fully explicit & Example~\ref{ex:rank1} \\
$S_2=G(1,1,2)$ & fully explicit & Example~\ref{ex:rank1} \\
$S_3=G(1,1,3)$ & fully explicit & Example~\ref{ex:S3} \\
$B_n=G(2,1,n)$ & $\fa$ generates $\Dc$ & Corollary~\ref{cor:36} \\
$G(m,1,n)$ ($e=1$) & $\fa$ generates $\Dc$ & Proposition~\ref{prop:36} \\
$D_3=G(2,2,3)$ & $\fa$ generates $\Dc$ & Lemma~\ref{lem:PsiTheta} (unconditional) \\
$D_n$, $n\geq4$ & open & needs $\fa^0$ with $\Theta,\Psi$ \\
$G(m,e,n)$, $e>1$, $\neq D_3$ & open & needs $\fa^0$ with $\Theta,\Psi$ \\
$F_4,H_3,H_4,E_6,E_7,E_8$ & $\fa$ generates $\Dc$, $\fa^0$ open & Remark~\ref{rem:exceptional} \\
\hline
\end{tabular}
\end{center}

\noindent Wherever it appears as a row's status, ``$\fa$ generates
$\Dc$'' is unconditional
(Proposition~\ref{prop:LS}); what remains open is the explicit
description of $\fa^0$ and $R_\fa$, the comparison
$R_\fa\leftrightarrow R_c$ of Remark~\ref{rem:37}, and, specifically
for the full lowest-weight equivalence of Corollary~\ref{cor:main},
the hypothesis of Lemma~\ref{lem:step0} that $\he$, up to an
additive scalar, lies in $\fa^0$, verified only in the rank-one case.

Building on the framework of Bøgvad and K\"allstr\"om~\cite{BK}, we
have shown that the decomposition of $B=\mathbb C[\h]$ as a
$\Dc$-module is governed by the Cartan-type ring $R_\fa$, via a
triangular structure on a Lie subalgebra $\fa\subset\Dc$
(Section~\ref{sec:triangular}) needed to circumvent the absence of an
a priori PBW factorisation for the spherical subalgebra; the precise
relationship between $R_\fa$ and the more naive $R_c$ is discussed in
Remark~\ref{rem:37}. This correspondence was made explicit for
$G(m,e,n)$ with $e=1$ and for the classical Coxeter group $B_n$; for
$D_n=G(2,2,n)$, generation was resolved only for $D_3$, unconditionally via
Lemma~\ref{lem:PsiTheta} (the exceptional isomorphism $D_3\cong A_3$
gives an independent, but conditional, route to the same conclusion),
with the general case
$n\geq4$ left open (see the table above). The unconditional part of
this contribution is methodological: adapting the lowest-weight
machinery of~\cite{BK} to a setting where the base algebra $\Dc$
itself has no a priori PBW structure (\S\ref{sec:triangular}) is what
makes the generation statements in the table possible at all, via
Levasseur--Stafford, for every family listed; the further step of
computing $\fa^0$ and $R_\fa$ explicitly, and with it the full
lowest-weight equivalence of Corollary~\ref{cor:main} without
appealing to Lemma~\ref{lem:step0}'s hypothesis, is carried out
unconditionally only in rank one, and is what remains open
for every other family treated.

Several questions remain open. First, the case $e>1$ of $G(m,e,n)$,
other than $D_3=G(2,2,3)$, requires a description of $\fa^0$ beyond
Proposition~\ref{prop:36}, since the basic invariants $f_i(x)$ no
longer freely generate $\mathrm{sym}(\h^*)^{G(m,e,n)}$; no exceptional
isomorphism is available to resolve this for general $D_n$ as it was
for $D_3$, nor for other $(m,e,3)$ with $e>1$. Short of an exceptional
isomorphism, one possible route
would be to enlarge $\fa$ from the outset to include $\Theta,\Psi$
among its generators (Proposition~\ref{prop:LS} places no restriction
on which homogeneous generators of $\mathrm{sym}(\h)^W$ and
$\mathrm{sym}(\h^*)^W$ are used to build $\fa$); whether the resulting
Lie algebra remains tractable enough to compute $\fa^0$ explicitly for
$D_n$, $n\geq4$, is not something we have investigated. Second, the compatibility of genericity of $c$ between $G$
and $\bar G$ for $e>1$ (Proposition~\ref{prop:36}(2)) has not been
established. Third, the injectivity of $R_\fa\to R_c$
(Remark~\ref{rem:37}) remains open beyond the rank-one case. Fourth,
the hypothesis of Lemma~\ref{lem:step0} that $\he$, up to an additive
scalar, lies in $\fa^0$ (needed for the full lowest-weight
equivalence of Corollary~\ref{cor:main}, though not for the weaker
generation statements in the table above) has likewise only been
verified in the rank-one case; a general argument (or a
counterexample) would be valuable. Finally,
it would be worthwhile to compare the resulting classification of
simple $\Dc$-modules with the classification, via category
$\mathcal O$~\cite{GGOR,Rouquier}, of simple modules over the full
Cherednik algebra $H_c$. The closest existing analogue is Gordon's
construction of baby Verma modules~\cite{Gordon}, quotients of
standard modules by their unique maximal submodule, though for the
degenerate parameter $t=0$ rather than the $t=1$ convention used
throughout this paper (a different algebra); making the relationship
between $eM(\tau)$ here and the $t=0$ baby Verma modules precise
would be worthwhile.

\section*{Acknowledgments}

We are deeply thankful to Professor Rikard Bøgvad for instructive
comments during the writing of this note.

\section*{Compliance with Ethical Standards}
\subsection*{Funding}
No funding was received to support this study 
\subsection*{Conflicts of Interest}
The authors declare that there are no conflicts of interest.
\subsection*{Author contribution declaration}  All the authors contributed equally in the conception and writing of this paper
\subsection*{Data Availability Statements}
No data were used to support this study

\end{document}